\documentclass{ifacconf}

\usepackage{graphicx}

\usepackage[labelformat=simple]{subcaption}

\usepackage{xcolor}
\usepackage{amsmath} 
\usepackage{mathrsfs}
\usepackage{amssymb}
\usepackage{mathtools}
\usepackage{mathdots}

\usepackage[shortlabels]{enumitem}
\usepackage{algorithm}
\usepackage[noend]{algpseudocode}
\usepackage{bm}
\usepackage{url}
\usepackage{optidef}
\usepackage{multirow}

\newtheorem{theorem}{\bf Theorem}[section]

\newtheorem{lemma}{\bf Lemma}[section]

\newtheorem{definition}{\bf Definition}[section]

\newtheorem{problem}{\bf Problem}[section]
\newtheorem{example}{\bf Example}[section]
\newtheorem{assumption}{\bf Assumption}[section]

\newcommand{\real}{{\mathbb{R}}}

\newcommand{\integersnonnegative}{\mathbb{Z}_{\geq 0}}
\newcommand{\integerspositive}{\mathbb{Z}_{> 0}}

\newcommand{\oprocendsymbol}{\hbox{$\bullet$}}
\newcommand{\oprocend}{\relax\ifmmode\else\unskip\hfill\fi\oprocendsymbol}
\newcommand{\proofendsymbol}{\hbox{$\blacksquare$}}
\newcommand{\proofend}{\relax\ifmmode\else\unskip\hfill\fi\proofendsymbol}

\newcommand{\defeq}{\vcentcolon=}

\newcommand{\rank}{{\textup{rank}}}

\newcommand{\Tf}{{T_{\textup{f}}}}

\newcommand{\hankel}[1]{\mathcal{H}_{#1}}
\usepackage{natbib}        % required for bibliography

\begin{document}
\begin{frontmatter}

\title{Robust Fundamental Lemma for Data-driven Control} 

\author[First]{Jeremy Coulson} 
\author[Second]{Henk van Waarde} 
\author[First]{Florian D\"orfler}

\address[First]{Automatic Control Laboratory, ETH Z\"urich, ZH 8006, Switzerland
(e-mail: \{jcoulson, dorfler\}@control.ee.ethz.ch).}
\address[Second]{Systems, Control and Optimization group, University of Groningen, 9747 AG Groningen, The Netherlands
(e-mail: h.j.van.waarde@rug.nl).}

\begin{abstract}                % Abstract of not more than 250 words.
The fundamental lemma by Willems and coauthors facilitates a parameterization of all trajectories of a linear time-invariant system in terms of a single, measured one. This result plays an important role in data-driven simulation and control. Under the hood, the fundamental lemma works by applying a persistently exciting input to the system. This ensures that the Hankel matrix of resulting input/output data has the ``right'' rank, meaning that its columns span the entire subspace of trajectories. However, such binary rank conditions are known to be fragile in the sense that a small additive noise could already cause the Hankel matrix to have full rank. Therefore, in this extended abstract we present a robust version of the fundamental lemma. The idea behind the approach is to guarantee certain lower bounds on the singular values of the data Hankel matrix, rather than mere rank conditions. This is achieved by designing the inputs of the experiment such that the minimum singular value of a deeper input Hankel matrix is sufficiently large. This inspires a new quantitative and robust notion of persistency of excitation. The relevance of the result for data-driven control will also be highlighted through comparing the predictive control performance for varying degrees of persistently exciting data. 
\end{abstract}


\end{frontmatter}

%===============================================================================
%===============================================================================
%===============================================================================

\section{Introduction}

The fundamental lemma from~\cite{JCW-PR-IM-BDM:05} is a powerful result that enables the characterization of the subspace of all possible trajectories of a linear time-invariant (LTI) system using raw time series data sorted into a Hankel matrix. The result has inspired many methods for data-driven analysis and control; see~\citep{IM-PR:08,HJVW-CDP-KMC-PT:20,HJVW:21,JC-JL-FD:18,CDP-PT:19,JB-JK-MM-FA:20}, and the survey by~\cite{IM-FD:21} and references therein.

At its core, the fundamental lemma requires applying an input sequence to the system that is persistently exciting of sufficient order such that the resulting input/output data Hankel matrix spans the entire subspace of possible trajectories of the system. In other words, if a deeper input data Hankel matrix is full row rank, the resulting input/output data Hankel matrix has low rank equal to the dimension of the admissible trajectory subspace. This input/output data Hankel matrix can then be used as a non-parametric system model. 

One major drawback of the fundamental lemma is that it only holds for noise-free data. Indeed, when the data are corrupted by noise, the input/output data Hankel matrix may no longer have low rank and thus will no longer describe the correct trajectory subspace leading to poor performance when used for data-driven analysis and control. This demonstrates that these binary rank conditions can no longer indicate suitable data when in the presence of noise. This motivates defining a new \emph{quantitative} notion of persistency of excitation that gives a measure of \emph{how} persistently exciting the inputs are. In adaptive control, such quantitative notions are studied~\citep[Remark 1, pg. 64]{KJA-BW:08-book}, but not in the context of the fundamental lemma.

In this extended abstract we propose a robust fundamental lemma that relies on a new \emph{quantitative} notion of persistency of excitation. The result informs the input selection such that the minimum singular value of the input/output data matrix is lower bounded by a user specified parameter. This results in a data matrix that is more robust to noise leading to better performance when used for data-driven analysis and control.

The rest of the extended abstract is organized as follows. We begin with notation. Section~\ref{sec:problem} formalizes the problem of interest. Section~\ref{sec:main} contains the main result whose relevance is illustrated through a data-driven control case study in Section~\ref{sec:examples}. We conclude in Section~\ref{sec:conclusion}.

\textit{Notation}: Let $m,n\in\integerspositive$. Given a matrix $M\in\real^{m\times n}$ and integer $i\in\integerspositive$ we denote the $i$-th singular value of $M$ by $\sigma_i(M)$ with the ordering $0\leq \sigma_1(M)\leq\sigma_2(M)\leq \cdots \leq \sigma_{\min\{m,n\}}$. 
When $m=n$ and $M=M^\top$, we denote the $i$-th eigenvalue of $M$ by $\lambda_i(M)$ with the ordering $\lambda_1(M)\leq \lambda_2(M)\leq\cdots \lambda_n(M)$. We use $\|M\|=\sigma_{\min\{m,n\}}$ to denote the spectral norm of a matrix $M$, $\|M\|_F$ the Frobenius norm, and $\|x\|$ to denote the 2-norm of a vector $x$. The Moore-Penrose pseudo inverse of $M$ is denoted by $M^\dagger$. Given $i,j,T\in\integersnonnegative$ with $i\leq j$ and a sequence $\{z(t)\}_{t=0}^{T-1}\subset\real^n$ define 
\[
z_{[i,j]}\defeq [z(i)^\top\;z(i+1)^\top\cdots z(j)^\top]^\top.
\]
Given an integer $k\in\integerspositive$ with $k\leq j-i+1$, define the \emph{Hankel matrix of depth $k$} associated with $z_{[i,j]}$ as
\[
\hankel{k}(z_{[i,j]}) \defeq 
\begin{bmatrix}
z(i) & z(i+1) & \cdots & z(j-k+1) \\
z(i+1) & z(i+2) & \cdots & z(j-k+2) \\
\vdots & \vdots & \ddots & \vdots \\
z(i+k-1) & z(i+k) & \cdots & z(j)
\end{bmatrix}.
\]

%===============================================================================
%===============================================================================
%===============================================================================
\section{Problem Statement}\label{sec:problem}
Consider the discrete-time LTI system
\begin{equation}\label{eq:sys}
x(t+1) = Ax(t)+Bu(t),
\end{equation}
where $x(t)\in\real^n$ is the state, and $u(t)\in\real^m$ is the control input. Throughout this extended abstract we will focus on a special case of full state measurements. We recall the definition of persistency of excitation.
\begin{definition}
Let $T\in \integerspositive$. The input sequence $u_{[0,T-1]}$ is called persistently exciting of order $k\in\integerspositive$ if $\hankel{k}(u_{[0,T-1]})$ has full row rank.
\end{definition}
We now state a version of the fundamental lemma~\citep{JCW-PR-IM-BDM:05} for the special case of input/state data.
\begin{theorem}\label{thm:fundamental}
Let $(A,B)$ be controllable and $T\in \integerspositive$. Let $(u_{[0,T-1]},x_{[0,T-1]})$ be an input/state trajectory of ~\eqref{eq:sys} such that $u_{[0,T-1]}$ is persistently exciting of order $n+1$. Then the data matrix
\begin{equation}\label{eq:HxHu}
\begin{bmatrix}\hankel{1}(x_{[0,T-1]})\\
\hankel{1}(u_{[0,T-1]})\end{bmatrix} = 
\begin{bmatrix} x(0) & x(1) & \cdots & x(T-1) \\
u(0) & u(1) & \cdots & u(T-1)
\end{bmatrix}
\end{equation}
has rank $n+m$.
\end{theorem}
As a result, the matrix~\eqref{eq:HxHu} along with $x(T)$ fully capture the behaviour of~\eqref{eq:sys} and can be used for data-driven analysis and control or system identification~\citep{IM-FD:21}.
However, the above rank condition on the data matrix~\eqref{eq:HxHu} is not always a valid indicator that the data is suitable for characterizing the behaviour of~\eqref{eq:sys}. In fact, when the data is corrupted by noise,~\eqref{eq:HxHu} may have full rank, but can lead to poor performance when used for data-driven analysis or control. We present a motivating example outlining precisely how we can improve on such rank conditions.
\begin{example}\label{ex:motivating}
Suppose we wish to identify matrices $A,B$ of the system 
\begin{equation}\label{eq:noisysys}
x(t+1) = Ax(t)+Bu(t)+w(t),
\end{equation}
where $w(t)\in\real^n$ is a noise term. Let $(u_{[0,T-1]},x_{[0,T-1]})$ be an input/state trajectory of the noisy system such that $u_{[0,T-1]}$ is persistently exciting of order $n+1$. To identify $A,B$, we consider the least squares problem
\[
\underset{A,B}{\min}\quad \left\| \hankel{1}(x_{[1,T]})-\begin{bmatrix} A & B \end{bmatrix} 
\begin{bmatrix}
\hankel{1}(x_{[0,T-1]})\\
\hankel{1}(u_{[0,T-1]})
\end{bmatrix}\right\|_F
\]
with solution given by
\begin{equation}\label{eq:ABestimate}
\begin{bmatrix} \hat{A} & \hat{B} \end{bmatrix} \defeq \hankel{1}(x_{[1,T]})\begin{bmatrix}
\hankel{1}(x_{[0,T-1]})\\
\hankel{1}(u_{[0,T-1]})
\end{bmatrix}^\dagger.
\end{equation}
The data satisfies
\[
\hankel{1}(x_{[1,T]}) = \begin{bmatrix} A & B \end{bmatrix} 
\begin{bmatrix}
\hankel{1}(x_{[0,T-1]})\\
\hankel{1}(u_{[0,T-1]})
\end{bmatrix}+
\hankel{1}(w_{[0,T-1]}).
\]
Hence, the error of our estimate is given by
\[
\begin{aligned}
\left\|\begin{bmatrix} \hat{A} & \hat{B} \end{bmatrix}-\begin{bmatrix} A & B \end{bmatrix}\right\|
&=\left\| \hankel{1}(w_{[0,T-1]})\begin{bmatrix}
\hankel{1}(x_{[0,T-1]})\\
\hankel{1}(u_{[0,T-1]})
\end{bmatrix}^\dagger \right\| \\
&\leq\frac{\sigma_n(\hankel{1}(w_{[0,T-1]})}{\sigma_1\left(\begin{bmatrix}
\hankel{1}(x_{[0,T-1]})\\
\hankel{1}(u_{[0,T-1]})
\end{bmatrix}\right)}.
\end{aligned}
\]
Note that the estimation error can be arbitrarily large, regardless of whether~\eqref{eq:HxHu} is full row rank. What is important instead is the smallest singular value of~\eqref{eq:HxHu}. Indeed, if we wish to identify $A,B$ to within some specified error, we require that the smallest singular value of the data matrix~\eqref{eq:HxHu} is bounded below by some threshold depending on the noise term $w$. \oprocend
\end{example}
This example motivates us to develop a robust fundamental lemma that moves away from rank conditions by defining a quantitative notion of persistency of excitation which guarantees a lower bound on the minimum singular value of the data matrix~\eqref{eq:HxHu}. More formally, the goal of this extended abstract is to solve the following problem.
\begin{problem}\label{prob:main}
Let $\delta>0$. Design an input sequence $u_{[0,T-1]}$ such that the resulting data matrix satisfies
\begin{equation}\label{eq:delta}
\sigma_1\left(
\begin{bmatrix}\hankel{1}(x_{[0,T-1]})\\
\hankel{1}(u_{[0,T-1]})\end{bmatrix} 
\right)
\geq \delta.
\end{equation}
\end{problem}

%===============================================================================
%===============================================================================
%===============================================================================

\section{Main Result}\label{sec:main}
The two main ingredients of the fundamental lemma are controllability and persistency of excitation. To establish a robust fundamental lemma, we must develop quantitative notions of these two main ingredients. We start with persistency of excitation.%
\begin{definition}\label{def:alphaPE}
Let $T\in \integerspositive$, $\alpha>0$. The input sequence $u_{[0,T-1]}$ is called $\alpha$-persistently exciting of order $k\in\integerspositive$ if $\sigma_1(\hankel{k}(u_{[0,T-1]}))\geq\alpha$.
\end{definition}
This is a natural generalization of persistency of excitation since for any $\alpha>0$ an $\alpha$-persistently exciting signal of order $k$ is necessarily persistently exciting of order $k$.

We now focus on the second main ingredient of the fundamental lemma: controllability.
Let
\begin{equation}\label{eq:z}
M\defeq \begin{bmatrix}
A & B & 0 & \cdots & 0 \\
0 & 0 & I_m & \cdots & 0 \\
\vdots & \vdots & \vdots & \ddots & \vdots \\
0 & 0 & 0 & \cdots & I_m \\
0 & 0 & 0 & \cdots & 0
\end{bmatrix},\quad z=\begin{bmatrix}\xi \\ \eta \\ 0_{nm}\end{bmatrix}
\end{equation}
where $\xi\in\real^n$, $\eta\in\real^m$, and $I_m$ and $0_{nm}$ denote the $m\times m$ identity matrix and the zero vector of length $nm$, respectively.
Define for any $z$ of the form~\eqref{eq:z} 
\begin{equation}\label{eq:theta}
\Theta_z \defeq \begin{bmatrix}
z & M^\top z & \cdots & (M^\top)^n z
\end{bmatrix}^\top.
\end{equation}
Note that $\Theta_z^\top$ can be viewed as an extended controllability matrix of the system $(M^\top,z)$. We make the following quantitative controllability assumption on this extended system.
\begin{assumption}\label{ass:rho}
Let $\rho>0$. For all $z$ of the form~\eqref{eq:z} with $\|z\|=1$, assume
\begin{equation}\label{eq:ass}
\sigma_1(\Theta_z)\geq \rho.
\end{equation}
\end{assumption}
This assumption can be thought of as ensuring that the finite horizon controllability Gramian $\Theta_z^\top\Theta_z=\sum_{j=0}^{n} (M^\top)^j z z^\top M^j$ of the extended system has lower bounded eigenvalues. This assumption is not restrictive since for any controllable $(A,B)$, there exists uniform $\rho_0>0$ that lower bounds $\sigma_1(\Theta_z)$ for all $z$, by the following lemma.
\begin{lemma}\label{lem:theta}
Let $(A,B)$ be controllable. Then 
\begin{enumerate}[(i)]
\item $\rank(\Theta_z)=n+1$, and \label{item:rank}
\item there exists $\rho_0>0$ such that $\sigma_1\left(\Theta_z\right)\geq\rho_0$, \label{item:rho0}
\end{enumerate}
for all $z$ be of the form~\eqref{eq:z} with $\left\|z\right\|=1$.
\end{lemma}
\begin{pf} We begin by proving~\ref{item:rank}. Assume on the contrary that the rows of $\Theta_z$ are linearly dependent. From the structure of $\Theta_z$, we must have that $\eta=0$ and thus $z=(\xi,0_{(n+1)m})$. Likewise, we must have that $\xi^\top B=\xi^\top AB=\cdots=\xi^\top A^{n-1}B=0$. By controllability of $(A,B)$, this implies that $\xi=0$, and hence $z=0$. This contradicts the fact that $\|z\|=1$, proving the claim. We now prove~\ref{item:rho0}. Assume on the contrary that $\forall \rho_0 >0$, $\exists \;z$ with $\|z\|=1$ such that $\sigma_1(\Theta_z)<\rho_0$. Let $\{\rho_j\}_{j=1}^\infty\subset \real$ be a sequence such that $\lim_{j\to\infty} \rho_j =0$. By assumption, for each $\rho_j$, there exists $z_j$ with $\|z_j\|=1$ and $\sigma_1(\Theta_{z_j})<\rho_j$. Since the set $\{z\mid  \|z\|=1\}$ is compact, $\{z_j\}_{j=1}^\infty$ has a convergent subsequence converging to some $\bar{z}$ with $\|\bar{z}\|=1$. Thus, $\sigma_1(\Theta_{\bar{z}})\leq0$. However, by part~\ref{item:rank}, $\Theta_{\bar{z}}$ has full row rank and hence, there exists $\bar{\rho}>0$ such that $\sigma_1(\Theta_{\bar{z}})>\bar{\rho}$ which is a contradiction, proving the result.
\proofend
\end{pf}

We now state the main theorem which solves Problem~\ref{prob:main}.
\begin{theorem}\label{thm:main}
Let $T\in \integerspositive$, $(A,B)$ be controllable, and $\delta>0$. Suppose that Assumption~\ref{ass:rho} holds and $(u_{[0,T-1]},x_{[0,T-1]})$ is an input/state trajectory of ~\eqref{eq:sys} such that $u_{[0,T-1]}$ is $\delta\frac{\sqrt{n+1}}{\rho}$-persistently exciting of order $n+1$. Then
\[
\sigma_1\left(\begin{bmatrix}\hankel{1}(x_{[0,T-1]}) \\ \hankel{1}(u_{[0,T-1]})\end{bmatrix}\right)\geq \delta.
\]
\end{theorem}
\begin{pf}
Denote the minimum singular value of the input/state data matrix~\eqref{eq:HxHu} by $\sigma\geq 0$ with corresponding left and right singular vectors $(\xi,\eta)\in\real^{n+m}$ and $v\in\real^T$, respectively. Then, $\left\|(\xi,\eta)\right\|=\|v\|=1$ and
\[
[\xi^\top\;\eta^\top]\begin{bmatrix}\hankel{1}(x_{[0,T-1]}) \\ \hankel{1}(u_{[0,T-1]})\end{bmatrix}=\sigma v^\top.
\]
Let $z=(\xi,\eta,0_{nm})$. By definition of $\Theta_z$ in~\eqref{eq:theta} and the dynamics~\eqref{eq:sys}, we have
\begin{equation}\label{eq:theta}
\Theta_z\begin{bmatrix}\hankel{1}(x_{[0,T-n-1]}) \\ \hankel{n+1}(u_{[0,T-1]})\end{bmatrix}
= \sigma \hankel{n+1}(v_{[0,T-1]}).
\end{equation}
By Cauchy's interlacing theorem~\cite[Corollary 3.1.3]{RAH-CRJ:94-book}, 
\begin{equation}\label{eq:ineq1}
\sigma_1\left(\hankel{n+1}(u_{[0,T-1]})\right)
\leq \sigma_{n+1}\left(\begin{bmatrix}\hankel{1}(x_{[0,T-n-1]}) \\ \hankel{n+1}(u_{[0,T-1]})\end{bmatrix}\right).
\end{equation}
By the Courant-Fischer-Weyl max-min principle~\cite[Theorem 3.1.2]{RAH-CRJ:94-book},
\[
\begin{aligned}
&\sigma_{n+1}\left(\begin{bmatrix}\hankel{1}(x_{[0,T-n-1]}) \\ \hankel{n+1}(u_{[0,T-1]})\end{bmatrix}\right) \\
&= \underset{\substack{\mathcal{U}\\ \dim(\mathcal{U})=n+1}}\min \;\underset{\substack{y\in\mathcal{U}\\ \|y\|=1}}\max \left\|\begin{bmatrix}\hankel{1}(x_{[0,T-n-1]}) \\ \hankel{n+1}(u_{[0,T-1]})\end{bmatrix}^\top y \right\|.
\end{aligned}
\]
By Lemma~\ref{lem:theta}, $\rank(\Theta_z)=n+1$, and hence
\begin{align}
&\sigma_{n+1}\left(\begin{bmatrix}\hankel{1}(x_{[0,T-n-1]}) \\ \hankel{n+1}(u_{[0,T-1]})\end{bmatrix}\right) \nonumber \\
&\leq \underset{\substack{y\in\textup{im}\Theta_z^\top \\ \|y\|=1}}\max \left\|\begin{bmatrix}\hankel{1}(x_{[0,T-n-1]}) \\ \hankel{n+1}(u_{[0,T-1]})\end{bmatrix}^\top y \right\| \nonumber\\
&=\underset{\substack{q \\ \|\Theta_z^\top q\|=1}}\max \left\|\sigma \hankel{n+1}(v_{[0,T-1]})^\top q \right\| \nonumber\\
&\leq \sigma \left\| \hankel{n+1}(v_{[0,T-1]})^\top \right\| \underset{\substack{q \\ \|\Theta_z^\top q\|=1}} \max \|q \| \nonumber \\
&\leq \sigma \sqrt{n+1} \underset{\substack{q \\ \|\Theta_z^\top q\|=1}} \max \|q \|  \label{eq:ineq2},
\end{align}
where the equality comes from~\eqref{eq:theta} and the last inequality holds because the rows of $\hankel{n+1}(v_{[0,T-1]})$ have norm at most $1$. Without loss of generality, write $q=\sum_{j=1}^{n+1} \alpha_j \mu_j$ where $\alpha_j\in\real$ and $\mu_j\in\real^{n+1}$ are orthonormal eigenvectors of $\Theta_z\Theta_z^\top$ corresponding to eigenvalues $\lambda_j(\Theta_z\Theta_z^\top)$. Then
\[
\underset{\substack{q \\ \|\Theta_z^\top q\|=1}} \max \|q \|^2 = \underset{\substack{\alpha_j \\ \sum_{j=1}^{n+1}\alpha_j^2\lambda_j(\Theta_z\Theta_z^\top)=1}} \max \sum_{j=1}^{n+1}\alpha_j^2.
\]
We see the maximum is achieved for $\alpha_1= \pm\frac{1}{\sqrt{\lambda_1(\Theta_z\Theta_z^\top)}}$, and $\alpha_j=0$ for $j\in\{2,\dots,n+1\}$. Hence,
\[
\underset{\|\Theta_z^\top q\|=1}\max \|q\|=\frac{1}{\sqrt{\lambda_1(\Theta_z\Theta_z^\top)}}.
\]
Combining the above with~\eqref{eq:ineq1} and~\eqref{eq:ineq2}, we obtain
\[
\begin{aligned}
&\sigma \geq \sigma_1\left(\hankel{n+1}(u_{[0,T-1]})\right)\sqrt{\frac{\lambda_1(\Theta_z\Theta_z^\top)}{n+1}}.
\end{aligned}
\]
Substituting $\sqrt{\lambda_1(\Theta_z\Theta_z^\top)} = \sigma_1(\Theta_z)\geq\rho$ by Assumption~\ref{ass:rho} and using the fact that $u_{[0,T-1]}$ is $\delta\frac{\sqrt{n+1}}{\rho}$-persistently exciting of order $n+1$ yields the result.
\proofend
\end{pf}
The theorem tells us how the inputs should be chosen such that \emph{for any} user defined parameter $\delta$, the smallest singular value of the data matrix~\eqref{eq:HxHu} is lower bounded by $\delta$. The degree of persistency of excitation needed depends on $\rho$ which is assumed to be a prior in our setting. However, Lemma~\ref{lem:theta} shows that, for any controllable system~\eqref{eq:sys}, there exists $\rho_0>0$ such that $\sigma_1(\Theta_z)\geq \rho_0$ for all $z$. To design an input sequence so that~\eqref{eq:delta} holds, we only require a lower bound on $\rho_0$ since any input that is $\delta\frac{\sqrt{n+1}}{\rho_0}$-persistently exciting of order $n+1$ yields the desired result.

%===============================================================================
%===============================================================================
%===============================================================================

\section{Numerical Example}\label{sec:examples}
\begin{figure}[h]
	\centering
		\includegraphics[width=\linewidth]{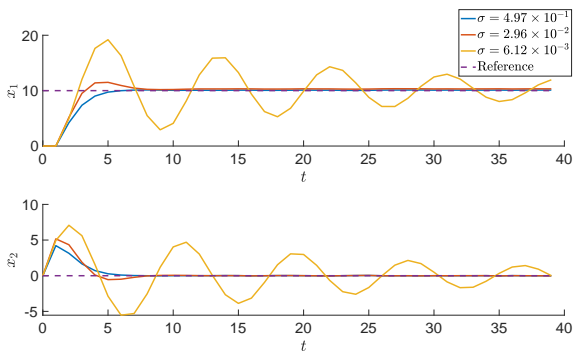}
	\caption{State trajectories $x(t)=(x_1(t),x_2(t))$ controlled by solving~\eqref{eq:deepc} in a receding horizon fashion for three data sets with varying degrees of persistency of excitation. The legend indicates the minimum singular value of the data matrices~\eqref{eq:HxHu} used for prediction in~\eqref{eq:deepc}.}
	\label{fig:trajectories}
\end{figure}

In this section we compare the performance of several data sets with varying degrees of persistency of excitation for a data-driven control task. Our hypothesis is that data sets whose inputs have a larger degree of persistency excitation will perform better when the data is corrupted by noise. Consider a controllable system~\eqref{eq:noisysys} with
\[
A = \begin{bmatrix} 1 & 1 \\ 0 & 1\end{bmatrix},\; B = \begin{bmatrix} 0 \\ 1 \end{bmatrix}.
\]
We generated 3 input data sequences of length $T=50$ with varying degrees of persistency of excitation. We denote the $i$-th input data sequence by $u^{(i)}_{[0,T-1]}$. The inputs for each input data sequence were chosen as $u^{(1)}(t) \sim \mathcal{N}(0,10^{-2}), u^{(2)}(t) =0.05u^{(1)}(t), u^{(3)}(t)=0.01u^{(1)}(t)$ for all $t\in\{0,\dots,T-1\}$. As a result, we obtained 3 input sequences that were $\alpha^{(i)}$-persistently exciting of order $n+1$, with $\alpha^{(1)}=0.48$, $\alpha^{(2)}=0.024$, $\alpha^{(3)}=0.0048$. The corresponding state data $x^{(i)}_{[0,T]}$ was generated by~\eqref{eq:noisysys} where the noise $w(t)\sim\mathcal{N}(0,10^{-4}I_n)$ was the same across all data sets. Using the data, we constructed 3 different data-driven one-step predictors as in~\eqref{eq:ABestimate}
\begin{equation}\label{eq:predictor}
x(t+1) = \hankel{1}\left(x^{(i)}_{[1,T]}\right)
\begin{bmatrix}
\hankel{1}\left(x^{(i)}_{[0,T-1]}\right)\\ 
\hankel{1}\left(u^{(i)}_{[0,T-1]}\right)
\end{bmatrix}^\dagger
\begin{bmatrix}
x(t) \\
u(t)
\end{bmatrix}.
\end{equation}
We then used the predictor equations~\eqref{eq:predictor} for a predictive control reference tracking task by solving the following optimization problem in a receding horizon fashion:
\begin{equation}\label{eq:deepc}
\begin{aligned}
\underset{x,u}{\min}\quad
&\sum_{k=0}^{\Tf-1} \|x_k-r\|^2+\|u_k\|^2 \\
\text{s.t.\quad}
& x_0 = x(t) \\
&x_{k+1} = \hankel{1}\left(x^{(i)}_{[1,T]}\right)
\begin{bmatrix}
\hankel{1}\left(x^{(i)}_{[0,T-1]}\right)\\ 
\hankel{1}\left(u^{(i)}_{[0,T-1]}\right)
\end{bmatrix}^{\dagger}
\begin{bmatrix}
x_k \\
u_k
\end{bmatrix}\\ 
& k\in\{0,\dots,\Tf-1\}
\end{aligned}
\end{equation}
where $x(t)$ denotes the current state at time $t$, $r\in\real^n$ is the reference, and $\Tf=10$ is the prediction horizon. System~\eqref{eq:noisysys} was simulated for each data set with noise $w(t)\sim\mathcal{N}(0,10^{-4}I)$ the same across all simulations. Figure~\ref{fig:trajectories} depicts the performance of the three data-driven one-step predictors for the reference tracking task.

For this example, Lemma~\ref{lem:theta} holds for $\rho_0=0.105$ and thus Assumption~\ref{ass:rho} holds for any $\rho\leq\rho_0$. Based on $\rho_0$ and the degree of persistency of excitation of each input sequence, we can compute the value $\delta^{(i)}$ for which~\eqref{eq:delta} holds for the $i$-th input/state data. Theorem~\ref{thm:main} guarantees that~\eqref{eq:delta} holds for each data set with $\delta^{(1)}=2.9\times10^{-2}$, $\delta^{(2)}=1.4\times 10^{-3}$, and $\delta^{(3)}=2.9\times 10^{-4}$. As we see from Figure~\ref{fig:trajectories}, the performance of the data-driven predictors decreases as the minimum singular value of the data matrix~\eqref{eq:HxHu} decreases. This can be attributed to the fact that the one-step predictor in~\eqref{eq:predictor} becomes more accurate as the smallest singular value of the data matrix decreases (as seen in Example~\ref{ex:motivating}), leading to better performance. 

%===============================================================================
%===============================================================================
%===============================================================================
\section{Conclusion}\label{sec:conclusion}
In this extended abstract, we defined a new quantitative notion of persistency of excitation by analyzing the minimum singular value of an input Hankel matrix. Based on this notion, we specified to what degree the inputs should be persistently exciting to ensure that the smallest singular value of the resulting input/state matrix is larger than a user defined threshold, thus generalizing the celebrated fundamental lemma. As a result, we are able to move away from classical rank conditions and give a quantitative notion of data suitability. A comparison of the control performance of several data sets with varying degrees of persistently exciting input data suggested that data being generated by inputs with a larger degree of persistency of excitation are better suited to control tasks. Future work includes extending these results to the general input/output case.
%===============================================================================
%===============================================================================
%===============================================================================

\bibliography{JC}             % bib file to produce the bibliography

\begin{thebibliography}{10}
\providecommand{\natexlab}[1]{#1}
\providecommand{\url}[1]{\texttt{#1}}
\providecommand{\urlprefix}{URL }
\expandafter\ifx\csname urlstyle\endcsname\relax
  \providecommand{\doi}[1]{doi:\discretionary{}{}{}#1}\else
  \providecommand{\doi}{doi:\discretionary{}{}{}\begingroup
  \urlstyle{rm}\Url}\fi

\bibitem[{{\AA}str{\"o}m and Wittenmark(2008)}]{KJA-BW:08-book}
{\AA}str{\"o}m, K.J. and Wittenmark, B. (2008).
\newblock \emph{Adaptive control}.
\newblock Dover Publications, 2nd edition.

\bibitem[{Berberich et~al.(2020)Berberich, K{\"o}hler, M{\"u}ller, and
  Allg{\"o}wer}]{JB-JK-MM-FA:20}
Berberich, J., K{\"o}hler, J., M{\"u}ller, M.A., and Allg{\"o}wer, F. (2020).
\newblock Data-driven model predictive control with stability and robustness
  guarantees.
\newblock \emph{IEEE Transactions on Automatic Control}, 66(4), 1702--1717.

\bibitem[{Coulson et~al.(2019)Coulson, Lygeros, and D{\"o}rfler}]{JC-JL-FD:18}
Coulson, J., Lygeros, J., and D{\"o}rfler, F. (2019).
\newblock Data-enabled predictive control: {I}n the shallows of the
  {D}ee{P}{C}.
\newblock In \emph{2019 18th European Control Conference (ECC)}, 307--312.

\bibitem[{De~Persis and Tesi(2019)}]{CDP-PT:19}
De~Persis, C. and Tesi, P. (2019).
\newblock Formulas for data-driven control: Stabilization, optimality and
  robustness.
\newblock \emph{IEEE Transactions on Automatic Control}.

\bibitem[{Horn and Johnson(1994)}]{RAH-CRJ:94-book}
Horn, R.A. and Johnson, C.R. (1994).
\newblock \emph{Topics in matrix analysis}.
\newblock Cambridge university press.

\bibitem[{Markovsky and D{\"o}rfler(2021)}]{IM-FD:21}
Markovsky, I. and D{\"o}rfler, F. (2021).
\newblock Behavioral systems theory in data-driven analysis, signal processing,
  and control.
\newblock \emph{Annual Reviews in Control}, 52, 42--64.

\bibitem[{Markovsky and Rapisarda(2008)}]{IM-PR:08}
Markovsky, I. and Rapisarda, P. (2008).
\newblock Data-driven simulation and control.
\newblock \emph{International Journal of Control}, 81(12), 1946--1959.

\bibitem[{van Waarde(2021)}]{HJVW:21}
van Waarde, H.J. (2021).
\newblock Beyond persistent excitation: Online experiment design for
  data-driven modeling and control.
\newblock \emph{IEEE Control Systems Letters}.

\bibitem[{van Waarde et~al.(2020)van Waarde, De~Persis, Camlibel, and
  Tesi}]{HJVW-CDP-KMC-PT:20}
van Waarde, H.J., De~Persis, C., Camlibel, M.K., and Tesi, P. (2020).
\newblock Willems' fundamental lemma for state-space systems and its extension
  to multiple datasets.
\newblock \emph{IEEE Control Systems Letters}, 4(3), 602--607.

\bibitem[{Willems et~al.(2005)Willems, Rapisarda, Markovsky, and
  De~Moor}]{JCW-PR-IM-BDM:05}
Willems, J.C., Rapisarda, P., Markovsky, I., and De~Moor, B.L. (2005).
\newblock A note on persistency of excitation.
\newblock \emph{Systems \& Control Letters}, 54(4), 325--329.

\end{thebibliography}
                                                     % with bibtex (preferred)
                                                  
\end{document}